\documentclass[aap,MSNbibl,seceqn,dvips]{arximspdf}
\usepackage{graphicx}

\doi{10.1214/09-AAP661}
\volume{20}
\issue{5}
\pubyear{2010}
\firstpage{1831}
\lastpage{1853}

\makeatletter
\DeclareMathAlphabet\mathcaligr{OMS}{cmsy}{m}{n}
\renewcommand{\mathcal}{\mathcaligr}
\renewcommand{\epsilon}{\varepsilon}
\newcommand{\ra}{\rightarrow}
\newcommand{\implies}{\Longrightarrow}
\newcommand{\mod}{\operatorname{mod}}

\newtheorem{theo}{Theorem}[section]
\newtheorem{pro}[theo]{Proposition}
\newtheorem{lem}[theo]{Lemma}
\newtheorem{cor}[theo]{Corollary}

\newproclaim{defin}[theo]{Definition}
\newproclaim{example}{Example}[section]
\newproclaim{prob}{Problem}
\newproclaim{con}{Conjecture}
\newproclaim{rem}[theo]{Remark}
\makeatother

\begin{document}
\begin{frontmatter}

\title{Time inhomogeneous Markov chains with~wave-like~behavior}
\runtitle{Time inhomogeneous Markov chains with wave-like behavior}

\begin{aug}
\author[A]{\fnms{L.} \snm{Saloff-Coste}\thanksref{t1}\ead[label=e1]{lsc@math.cornell.edu}}
\and
\author[B]{\fnms{J.} \snm{Z\'u\~niga}\corref{}\thanksref{t2}\ead[label=e2]{jzuniga@math.stanford.edu}}
\runauthor{L. Saloff-Coste and J. Z\'u\~niga}
\affiliation{Cornell University and Stanford University}
\address[A]{Department of Mathematics\\
Cornell University\\
Malott Hall\\
Ithaca, New York 14853\\
USA\\
\printead{e1}} 
\address[B]{Department of Mathematics\\
Stanford University\\
Building 380\\
Stanford, California 94305\\
USA\\
\printead{e2}}
\end{aug}
\pdfauthor{L. Saloff-Coste and J. Zuniga}
\thankstext{t1}{Supported in part by NSF
Grant DMS-06-03886.}
\thankstext{t2}{Supported in part by NSF
Grants DMS-06-03886, DMS-03-06194 and DMS-08-03018.}
\received{\smonth{1} \syear{2009}}
\revised{\smonth{7} \syear{2009}}

\begin{abstract}
Starting from a given Markov kernel on a finite set $V$
and a bijection $g$ of $V$,
we construct and study  a time inhomogeneous Markov chain  whose kernel at time $n$
is obtained from $K$ by transport of $g^{n-1}$.
We show that this construction leads to interesting examples, and we obtain
quantitative results for some of these examples.
\end{abstract}


\end{frontmatter}

\section{Introduction}  In \cite{SZ,SZ3,SZ4}, we considered the
problem of obtaining quantitative results describing the ergodic behavior
of time inhomogeneous  finite Markov chains. In general, a time inhomogeneous Markov
chain, say on a finite set $V$,  is described by a sequence of Markov kernels
$(K_i)_1^\infty$. At time $n$,
the distribution of the chain started at $x$ is denoted by
$K_{0,n}(x,\cdot)$. More generally, for $n\le m$, we define $K_{n,m}$ inductively by
$K_{n,n}=I$ (the identity matrix) and
\[
K_{n,m}(x,y)=\sum_{z}K_{n,m-1}(x,z)K_m(z,y), \qquad  x,y\in V .
\]

If each $K_i$ is irreducible and aperiodic, one expects that, in many cases, the
Markov chain driven by this sequence will have the property that
\[
\forall x,y \qquad
\|K_{0,n}(x,\cdot)-K_{0,n}(y,\cdot)\|_{\mathrm{TV}} \ra 0\qquad \mbox{as } n\ra \infty.
\]
We call this property \textit{total variation merging} and say that the
chain driven by the sequence
$(K_i)_1^\infty$ is merging.  Note that, in general, $K_{0,n}(x,\cdot)$ does not tend to
a limiting distribution. However, when merging occurs, the chain does forget where
it started: asymptotically, the distribution sequence evolves in time following a
well-defined pattern which is  independent of the starting distribution.

In this paper, we will mostly discuss a stronger notion which we call
relative-sup merging.  By definition, the sequence $(K_i)_1^\infty$
is merging in relative-sup if
\[
\max_{x,y,z\in V} \biggl\{ \bigg|\frac{K_{0,n}(x,z)}
{K_{0,n}(y,z)}-1\bigg | \biggr\}\ra 0 \qquad\mbox{as } n\ra \infty.
\]

In general, the relative-sup distance between two measures $\mu$ and $\nu$
(on a finite or countable state space) is defined by
(note the asymmetry)
\[
\max_{ x\in V} \biggl\{ \bigg|\frac{\mu(x)}{\nu(x)}-1 \bigg| \biggr\}.
\]
In particular, for a time inhomogeneous chain driven by a sequence
$(K_i)_1^\infty$ of Markov kernels, we will
consider quantities such as
\[
\max_{x,z\in V} \biggl\{ \bigg|\frac{K_{0,n}(x,z)}{\mu_n(z)}-1\bigg |
\biggr\},
\]
where $\mu_n=\mu_0K_{0,n}$ for some starting measure $\mu_0$.
For $\epsilon>0$, we also define  the
$\epsilon$ relative-sup merging time $T_\infty(\epsilon)$ by
\[
T_\infty(\epsilon)=\min \biggl \{n \dvtx\max_{x,y,z\in V}
\biggl \{ \bigg|\frac{K_{0,n}(x,z)}{K_{0,n}(y,z)}-1\bigg | \biggr\}<\epsilon \biggr\}.
 \]
See \cite{SZ3} for more details.

Background and general results concerning time inhomogeneous Markov chains
are described in \cite{Io,Pa,Sen} where further
references can be found. It turns out that the study of merging
is difficult, both at the qualitative and the quantitative level,
except in the special but interesting case when all the kernels in the sequence $(K_i)_1^\infty$ share
the same stationary probability measure. See, for example, \cite{DR,Ga,MPS,SZ}.
Only a small set of examples have been treated in the literature
mostly because proving  anything about concrete time inhomogeneous Markov chains is
difficult.

This  paper describes a  special class of examples whose structure is, in itself,
quite interesting and for which some results can be obtained.
The set up is as follows. On a finite or countable set $V$,  we are given a
Markov kernel $K$ and a bijection $g \dvtx V\ra V$. We then consider the time
inhomogeneous Markov chain driven by the sequence of the kernels
\[
K_i(x,y)=K(g^{i-1}x,g^{i-1}y), \qquad  x,y\in V,   i=1,2,\ldots .
\]
The problem is to study this time inhomogeneous chain and its merging properties. As
we shall see, this covers some interesting examples and leads
to interesting results as well as difficult open problems.

The examples discussed in this paper can serve to illustrate the techniques
developed in \cite{SZ3,SZ4}. In particular, we will make use of the following basic
singular value technique. See \cite{DLM} and  Theorem~3.2 of  \cite{SZ3}.

\begin{theo}\label{thm-sv-merge}
Given a sequence of Markov kernels $K_i$, $i=1,2,\ldots,$ on a set $V$ and a
positive probability measure $\mu_0$, set $\mu_n=\mu_0K_{0,n}$ and let
$\sigma_1(i)$ be the second largest singular value of the operator
$K_i \dvtx \ell^2(\mu_i)\ra \ell^2(\mu_{i-1})$. Then
\[
 \bigg|\frac{K_{0,n}(x,z)}{\mu_n(z)}-1 \bigg|\le
 \biggl(\frac{1}{\mu_0(x)}-1 \biggr)^{1/2} \biggl(\frac{1}{\mu_n(z)}-1 \biggr)^{1/2}
\prod_1^n \sigma_1(i).
\]\end{theo}

This good-looking result is deceptive because, unless one can get some
control on the sequence of measures $\mu_n$, it is essentially useless.
Note in particular that $\sigma_1(n)$ depends very much on
$\mu_{n-1}$ and $\mu_n$.

\section{Stability}
\label{sec-stab}

It is well established that the stationary distribution of an irreducible aperiodic
time homogeneous Markov chain plays a crucial part in the
analysis of the ergodic properties of the chain. Not much can be said unless one can
get some control on the stationary distribution. Moreover, unless the chain is
reversible or some algebraic miracle occurs, the computation of the stationary
measure is a difficult problem.

The situation for time inhomogeneous Markov chains is much worse. In order to
understand how the chain behaves when started from an arbitrary distribution, it is
crucial to find (at least) \textit{one} initial distribution $\mu_0$ such that sequence
of probability measures  $\mu_n=\mu_0K_{0,n}$ is somewhat well behaved. The ideal
situation is when there is a $\pi$ such $\pi K_{0,n}=\pi$. This occurs if an only if
all $K_i$ admit the same invariant
measure $\pi$, a rather fortunate but rare circumstance.
The next definition, taken from \cite{SZ3},
introduces a property that is an obvious weakening
of the existence of a common invariant measure.

\begin{defin} Fix $c\ge 1$.
A sequence of Markov kernels $(K_n)_1^{\infty}$ on a finite set $V$
is $c$-stable if there exists  a measure $\mu_0$  such that
\begin{equation}\label{c-stab-seq}
\forall n\ge 0, \  x\in V \qquad  c^{-1}\leq \frac{\mu_n(x)}{\mu_0(x)}\leq
c,
\end{equation}
where $\mu_n=\mu_0K_{0,n}$. If this holds, we say that
$(K_n)_1^{\infty}$
is $c$-stable with respect to the measure $\mu_0$.
\end{defin}

We refer the reader to \cite{SZ3,SZ4}, for examples, and results involving
$c$-stability. The idea behind this definition is that, if a sequence is $c$-stable
with respect to a probability measure $\mu_0$, then one can
study the merging of this sequence \textit{more or less} as one would study the
ergodicity of a time homogeneous chain with invariant measure $\mu_0$.
Why this is true is not obvious and the required technical details  are quite
intricate. Precise results in this direction are described in \cite{SZ3,SZ4}. We
think that $c$-stability is an interesting property in itself and that it deserves
some attention. Note also that, even for a fixed sequence $(K_i)_1^\infty$ on a
fixed finite state space, $c$-stability is a nontrivial property. The case of the
two point space is treated  in \cite{SZ3}.

A special case of interest to us here is when the time inhomogeneous Markov chain is
driven  by a sequence
$(K_i)_1^\infty$ that is periodic in the sense that there is an integer $k$ such that
\[
\forall i \qquad   K_{i+k}=K_i.
\]
In such case, there is an obvious candidate for a ``good'' starting distribution
$\mu_0$, namely,
the invariant measure  $\pi$ of $K_1\cdots K_k=K_{0,k}$. Indeed, if we pick
$\mu_0=\pi$ then
the sequence $\mu_n=\mu_0 K_{0,n}$
is also periodic of period $k$.  If we can compute~$\pi$, this might
allow us to investigate the property of the sequence $\mu_n$ including
$c$-stability. Note however that in many examples of interest,
the period $k$ will grow with the size of the state space $V$ so that,
even in that case, investigating $c$-stability in a meaningful way is difficult.

An example of this type is cyclic to random transpositions. On $V=S_n$, the
symmetric group, let $Q_i$ be the Markov kernel $Q_i(x,y)=1/n$ if $y=x$ or if
$y=x(i,j)$ for some $j\neq i$ and $Q_i(x,y)=0$ otherwise.
Here $(i,j)$ stands for the corresponding  transposition. This kernel
corresponds to ``transpose the card in position $i$ with the card in a uniformly
chosen position.'' The cyclic-to-random transposition chain is driven by the
sequence of kernels $(K_i)_1^\infty$ with $K_{i}=Q_{i \mod n}$ (by definition,
$Q_0=Q_n$). See \cite{Ga,MPS,SZ}. Of course, in this example, the uniform measure is
invariant for all $Q_i$. Other examples of periodic time inhomogeneous chains
are discussed in \cite{DR}.

\section{Periodic waves}
\label{sec-wave}

We now describe in detail the construction outlined in the introduction.
This construction is of a rather general nature and
produces periodic time inhomogeneous Markov chains that reduce, in a sense,
to  time homogeneous chains.

Let $K$ be a Markov kernel on a finite  state space $V$, and let
$g \dvtx V\ra V$, $x\mapsto g(x)=gx$ be a bijection.  The order of the map $g$ is
\[
k=\min\{n\in\mathbb{N} \dvtx \forall x\in V  g^nx=x\},\qquad   g^n=g\circ g\circ \cdots\circ g.
\]
For all $x,y\in V$, set
\begin{equation}\label{g-i-kernel}
K_i(x,y)=K(g^{i-1}x,g^{i-1}y)
\end{equation}
so that $K=K_1$. Consider the inhomogeneous Markov chain driven by the
sequence $(K_i)_1^{\infty}$ defined above.
It is easy to see that  all  $K_i$ are irreducible aperiodic kernels
if and only if $K$ is. Moreover, if $K$ has stationary distribution $\pi$
then $K_i$ has stationary distribution $\pi_i$ where $\pi_i(x)=\pi(g^{i-1}x)$.
Obviously, the sequence $(K_i)_1^\infty$ is periodic of period $k$.
Examples are discussed below after we discuss some general properties of these
chains.  Given this definition, the obvious question we face is the following:  How
are the (quantitative) merging properties of the chain driven by $(K_i)_1^\infty$
related to the (quantitative) ergodic properties of the
chain driven by $K$?

\begin{pro}\label{pro-kernel-shift}
Set
\begin{equation}\label{kernel-shift}
\widetilde{K}(x,y)=K(x,g^{-1}y),
\end{equation}
where $g^{-1} \dvtx V\ra V$ is the inverse of the map $g$. Then
$K_{0,n}$ is given by
\[
K_{0,n}(x,y) =\widetilde{K}^n(x,g^ny).
\]
\end{pro}
\begin{pf} We proceed by induction. For $n=1$ the result holds by definition.
Assume that $\widetilde{K}^n(x,y)=K_{0,n}(x,g^{-n}y)$. Then we have
\begin{eqnarray*}
\widetilde{K}^{n+1}(x,y)
&=&\sum_{z\in V}\widetilde{K}^n(x,z)\widetilde{K}(z,y)\\
&=&\sum_{z\in V}K_{0,n}(x,g^{-n}z)K_{n+1}(g^{-n}z,g^{-n-1}y)\\
&=&K_{0,n+1}\bigl(x,g^{-(n+1)}y\bigr).
\end{eqnarray*}
This gives the desired result. \end{pf}

\begin{figure}[b]

\includegraphics{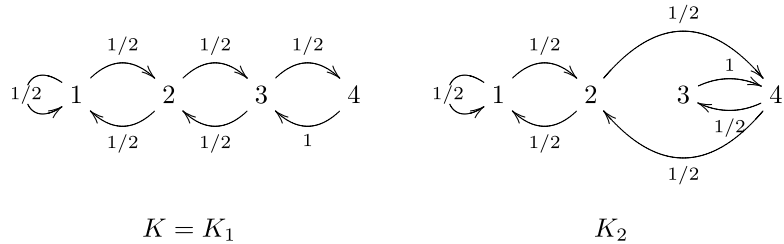}

\caption{Graph structure for kernels $K$ and $K_2$.\label{fig-kernels-4pt}}
\end{figure}

\begin{cor} \label{cor-irr}
The kernel $\widetilde{K}$ is irreducible aperiodic
if and only if there exists an integer $n_0>0$ such that for all $x,y\in V$,
$K_{0,n_0}(x,y)>0$.
\end{cor}

The following examples illustrate some of the subtleties of this construction.

\begin{example} Let $K$ be irreducible, periodic of period $k$,
with periodicity classes $C_0,\ldots,C_{k-1}$ so that $K(x,y)>0$ if and only if
$x\in C_i$ and $y\in C_{i+1\mod k}$. Assume that $|C_0|=\cdots=|C_{k-1}|$,
that is, all the periodicity classes have the same cardinality.
Let $g \dvtx V\ra V$ be a bijection such that $g(C_i)=C_{i-1 \mod k}$.
Let $K_i(x,y)=K(g^{i-1}x,g^{i-1}y)$,
$\widetilde{K}(x,y)=K(x,g^{-1}y)$ as above.
It is clear that $\widetilde{K}(x,y)>0$ if and only if
 $x,y$ are in the same class $C_i$ for some $i$. That is,
$\widetilde{K}$~is not irreducible. One the other hand,
for any $x,y$ there exists $n=n(x,y)$ such that $K_{0,n}(x,y)>0$.
\end{example}

\begin{example}\label{ex-4points}
On $V=\{1,2,3,4\}$, consider the irreducible aperiodic reversible kernel
$K$ given by $K(1,1)=K(1,2)=K(2,1)=K(2,3)=K(3,2)=K(3,4)=1/2$, $K(4,3)=1$
and $K(x,y)=0$, otherwise. Let $g$ be the map that transposes $3$ and $4$.
Then $K_2(1,1)=K_2(1,2)= K(2,1)=K_2(2,4)=K_2(4,2)=K_2(4,3)=1/2$,
$K_2(3,4)=1$ and $K_2(x,y)=0$, otherwise.  The graph structure for kernels
$K$ and $K_2$ is illustrated in Figure \ref{fig-kernels-4pt}.
It follows that
\[
K_{0,2n}(4,4)=1, \qquad   K_{0,2n+1}(4,3)=1.
\]
This shows that the property that $K$ is irreducible aperiodic does
not imply that for each $x,y$ there is an $n=n(x,y)$ such that
$K_{0,n}(x,y)>0$. Further,
$\widetilde{K}(1,1)=\widetilde{K}(1,2)=\widetilde{K}(2,1)=
\widetilde{K}(2,4)=\widetilde{K}(3,2)
=\widetilde{K}(3,3)=
1/2$, $\widetilde{K}(4,4)=1$. Hence,
$\widetilde{K}$~is not irreducible and has a unique absorbing state, namely,
the point $4$ as illustrated by Figure \ref{fig-tildeK}.

\begin{figure}

\includegraphics{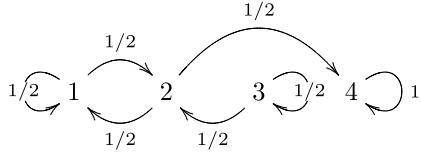}

\caption{Graph structure for $\widetilde{K}$.}\label{fig-tildeK}
\end{figure}

This implies that  the sequence $K_1,K_2,K_1,K_2,\ldots$
is merging in total variation, that is,
$K_{0,n}(x,z)-K_{0,n}(y,z)\ra 0$ for any $x,y,z$. Note that for $z\neq 4$,
we  have $K_{0,2n}(x,z)\ra 0$ for any $x$. However, this same sequence is not
merging in relative-sup distance. Indeed,
\[
T_\infty(\epsilon)=\min \biggl\{n \dvtx \max_{x,y,z} \biggl\{
 \bigg|\frac{K_{0,n}(x,z)}{K_{0,n}(y,z)}-1 \bigg| \biggr\}<\epsilon \biggr\}=\infty
 \]
since $K_{0,2n}(4,1)=0$ and $K_{0,2n}(1,1)>0$.

This gives an example of a pair $K_1,K_2$ of reversible,  irreducible and
aperiodic Markov kernels such the sequence $K_1,K_2,K_1,K_2,\ldots$ is not
merging in relative-sup distance.
\end{example}

\begin{example}
On the symmetric group $S_n$, set  $\sigma$ and $\sigma'$ to be the cycles
$\sigma=(n,n-1,\ldots,1)$ and $\sigma'=(n-1,n-2,\ldots,1)$ and $a$ to be the permutation
defined by $a(i)=n-i+1$.
In terms of a deck of $n$ cards, $\sigma$ takes the top card to the bottom,
$\sigma'$ takes the top card to the second to last position whereas $a$
reverses the order of the deck. Consider the kernel
$K(x,y)=1/2$ if  $x^{-1}y\in \{\sigma,\sigma'\}$ and $0$ otherwise,
and the bijection $g(x)= axa^{-1}$, which is of order $2$.
Observe that~$K$ is irreducible and aperiodic.
Note that $g (\sigma) = \sigma^{-1}$ (take the bottom card and put it on top)
and $g (\sigma')=(2,3,\ldots,n)$
(take the bottom card and put it in second position).
From this it follows that
\begin{eqnarray*}
K_{0,2}(x,y) &= &\sum_z K(x,z)K(g(z),g(y))\\
&=&
\cases{ 1/4,
&\quad  if   $x^{-1}y\in\{ e, (1,2), (1,n), (1,n,2)\}$,\cr
0, &\quad  otherwise.}
\end{eqnarray*}
This shows that, for all $n$,  $K_{0,2n}(e,x)=0$ unless
$x\in B=\{e, (1,2), (1,n), (1,n,\break2)\}$, and
$K_{0,2n+1}(e,x)=0$ unless
$x\in \sigma B\cup \sigma' B$. We note that describing $\widetilde{K}$
is difficult.
\end{example}

\begin{pro}
Let $\widetilde{\pi}$ be an invariant measure for $\widetilde{K}$. Set
\[
\forall x\in V,  i=1,2,\ldots \qquad   \mu_i(x)=\widetilde{\pi}(g^{i}x).
\]
Then $\mu_{i-1}K_i=\mu_{i}$.
\end{pro}

\begin{pf} Indeed, we have
\begin{eqnarray*}
\mu_{i-1}K_i(x)
&=&\sum_{z\in V}\mu_{i-1}(z)K_i(z,x)
=\sum_{z\in V}\widetilde{\pi}(g^{i-1}z)K_1(g^{i-1}z,g^{i-1}x)\\
&=&\sum_{z\in V}\widetilde{\pi}(g^{i-1}z)\widetilde K(g^{i-1}z,g^{i}x)
=\widetilde{\pi}(g^ix)=\mu_{i}(x).
\end{eqnarray*}
\upqed
\end{pf}
The ``wave'' appearing in the title of this paper corresponds to the
distribution~$\widetilde{\pi}$. The time inhomogeneous chain driven by the sequence
$(K_i)_1^\infty$  produces the wave~$\widetilde{\pi}$, moving around in a
periodic fashion under the action of the bijection~$g$ on the set $V$.
Despite the similarity in names, we do not claim any connection of
this paper with the subject of traveling waves.

\begin{cor}\label{mv-stab}
Assume that $\widetilde{K}$ admits a positive invariant measure $\widetilde{\pi}$.
Then the sequence $(K_n)_1^{\infty}$ is
$c$-stable with respect to the measure $\mu_0=\widetilde{\pi}$
with
\[
c=\max_{x,i}\{\widetilde{\pi}(g^ix)/\widetilde{\pi}(x)\}.
\]
\end{cor}

The next proposition discusses the singular value decompositions
of various operators appearing in this construction.
The proof is by inspection. We use the following notation.
We assume that $\widetilde{\pi}$ is an invariant measure for $\widetilde{K}$
and that $\widetilde{\pi}(x)>0$ for all $x \in V$.
Let $\widetilde{\sigma}_j$, $j=0,\ldots,|V|-1,$ be the singular values of
$\widetilde{K} \dvtx\ell^2(\widetilde{\pi})\ra \ell^2(\widetilde{\pi})$ in nonincreasing
order,
and let $(\widetilde{\phi}_j)_0^{|V|-1}$, $(\widetilde{\psi}_j)_0^{|V|-1}$,  be
orthonormal bases
of $\ell^2(\widetilde{\pi})$ such that
$\widetilde{K}\widetilde{\phi}_j=\widetilde{\sigma}_j \widetilde{\psi}_j$ (with
$\widetilde{\sigma}_0=1, \widetilde{\phi}_0=\widetilde{\psi}_0\equiv 1$).
We refer the reader to \cite{SZ3} for a detailed discussion.
The orthonormal  bases $(\widetilde{\phi}_j)_0^{|V|-1}$,
$(\widetilde{\psi}_j)_0^{|V|-1}$ are, respectively,  eigenbases  for
$K^*K$ and $KK^*$.

\begin{pro}\label{pro-sd-wave}
For any $i=1\in\{1,\ldots\}$,
$ \phi_j^i(x)= \widetilde{\phi}_j(g^{i}x)$, $j=0,\ldots,\break|V|-1$,
and $\psi_j^i(x)= \widetilde{\psi}_j(g^{i-1}x)$, $j=0,\ldots,|V|-1$,
are orthonormal bases of $\ell^2(\mu_{i})$ and
$\ell^2(\mu_{i-1})$, respectively,
which provide a singular value decomposition of $K_i \dvtx\ell^2(\mu_i)\ra
\ell^2(\mu_{i-1})$ in the sense that
$K_i\phi^i_j=\widetilde{\sigma}_j\psi_j^i$.
In particular, the singular values $\sigma_j(K_i,\mu_{i-1})$
of $K_i \dvtx\ell^2(\mu_i)\ra
\ell^2(\mu_{i-1})$  are given by
$\sigma_j(K_i,\mu_{i-1})
=\widetilde{\sigma}_j$, $j=0,\ldots,|V|-1$.

If $\widetilde{\alpha}$ is an eigenvalue of $\widetilde{K}$ with eigenfunction
$\widetilde{\omega}$ and $k$ is the order of $g$ then
$\widetilde{\alpha}^k$ is an eigenvalue of $K_1\cdots K_k$ with the same eigenfunction.
\end{pro}

This proposition illustrates clearly the difficulties that appear in relating the
ergodic properties of the kernel $K$ (that serves as the basic ingredient of this
construction) to
the merging properties of the sequence $(K_i)_1^\infty$.
Indeed, it is rather unclear how the
ergodic properties of $K$ and the properties of its stationary measure $\pi$ relate
to $(\widetilde{K},\widetilde{\pi})$.

In the following two examples,
$\pi=\widetilde{\pi}$ is the uniform measure on $V$. Even in these cases, the
above construction is quite interesting and nontrivial. Examples with
$\pi\neq \widetilde{\pi}$ will be discussed in the next two sections.

\begin{example}[(Cycling for binary vectors)] \label{ex-cyc-bv}
In this example, the kernel $K$ is not irreducible.
Take $V=\{0,1\}^N$ with  $\pi$ being the uniform distribution on $V$. Let~$e_i$ be
the binary vector with a unique $1$ in position $i$. Let $K(x,y)=0$ except if $y=x$
or $y=x+e_1$ in which case $K(x,y)=1/2$ ($K$ randomizes the first binary entry of
$x$). Let $gx=(x_2,\ldots,x_N,x_1)$ if $x=(x_1,\ldots,x_N)$ (shift to the left).
Using the definition, one checks that $K_i$ is the Markov kernel that
randomizes the $i$th coordinate. Hence, $K_1\cdots K_N=\pi$
(after $N$ steps, we have a binary vector picked uniformly at random).

The kernel $\widetilde{K}$ corresponds to randomizing the first entry and
shifting  left. Its invariant measure  $\widetilde{\pi}$ is  uniform. One recovers
immediately the fact that the uniform distribution is reached after exactly $N$
steps. The singular values ($=$eigenvalues) of~$K$ (which is reversible) are $1$
(multiplicity $2^N-1$) and $0$ (multiplicity $1$). The kernel $\widetilde{K}$ has the
property that $\widetilde{K}^*\widetilde{K}=K=K^2$ so that it has the same singular
values. The operator $\widetilde{K}$ has two eigenvalues, $0$ and $1$, and is not
diagonalizable, but $\widetilde{K}-\pi$
is nilpotent since $(\widetilde{K}-\pi)^N=0$.
\end{example}

\begin{example}[(Cyclic-to-random transposition)] See, for example, \cite{MPS,SZ}.
On the symmetric group $S_n$, let $K(x,y)=1/n$ if $y=x(1,j)$,
$j=1,2,\ldots,n$, and $K(x,y)=0$ otherwise (this is called
``transpose top with random'').
Let $\sigma$ be the cycle $(1,2,\ldots,n)$ and
$g \dvtx S_n\ra S_n$, $x\mapsto g(x)= \sigma x \sigma^{-1}$. Observe that
$g^i((1,j))= (i,j+i \mod n)$ so that $K_i$ is
``transpose $i$ with random.'' Hence, we recover the cyclic-to-random
transposition chain.

Because $\widetilde{\pi}=\pi$ in this case,
it follows that the singular values of $\widetilde{K}$ are equal to the
singular values of $K$ which can be computed by using the
representation theory of $S_n$.
Note that, as $K$ is reversible, the singular values of $K$
are the square roots of the square of its eigenvalues, that is,
the absolute value of
the eigenvalues.  In particular,
$\widetilde{\sigma}_1=1-1/n$ and thus $\sigma_1(K_i,\pi)=\widetilde{\sigma}_1
=1-1/n$ for all $i$ (see \cite{Dia2,FOW,SZ,SZ1}).
The eigenvalues of $\widetilde{K}$ are rather
mysterious, and it is not clear that $\widetilde{K}$ is diagonalizable. See
\cite{MPS} where the eigenvalues of $K_1\cdots K_n$
(hence, indirectly, the eigenvalues of $\widetilde{K}$) are investigated
and used to obtain a very interesting lower bound on the mixing time of
cyclic to random transposition.
\end{example}

Propositions \ref{pro-kernel-shift} and \ref{pro-sd-wave} reduce
the study of the merging\vspace*{2pt} of the sequence $(K_i)_1^\infty$  to
the study of the ergodicity of the  time homogeneous Markov chain driven by
$\widetilde{K}$.  More precisely, we have the following result.

\begin{theo} \label{th-wave}
Fix $V,K,g,\widetilde{K}$ and $(K_i)_1^\infty$ as above.
\begin{longlist}[(2)]
\item[(1)] The sequence $(K_i)_1^\infty$ is merging in relative-sup
if and only if the kernel $\widetilde{K}$ is irreducible and aperiodic.
\item[(2)] If  $\widetilde{K}$ is irreducible and aperiodic,
let $\widetilde{\pi}$ be its unique  invariant probability measure and  set
$\mu_i(x)=\widetilde{\pi}(g^{i}x)$, $x\in V$. Then
\[
 \bigg|\frac{K_{0,n}(x,z)}{\mu_n(z)}-1 \bigg|\le
 \biggl(\frac{1}{\widetilde{\pi}(x)}-1 \biggr)^{1/2}
 \biggl(\frac{1}{\widetilde{\pi}(g^nz)}-1 \biggr)^{1/2}
\widetilde{\sigma}_1 ^n,
\]
where $\widetilde{\sigma}_1$ is the second largest singular value of $\widetilde{K}$
acting on $\ell^2(\widetilde{\pi})$.
\end{longlist}
\end{theo}

\begin{pf} Use Propositions \ref{pro-kernel-shift} and \ref{pro-sd-wave}.
To obtain the last inequality, use Theorem~\ref{thm-sv-merge}. Theorem~3.2
of \cite{SZ3} also yields additional inequality for the chi-square distance between
$K_{0,n}(x,\cdot)$ and $\mu_n$.
\end{pf}

\begin{rem}
Example \ref{ex-4points} gives an example where total variation merging
occurs, but $\widetilde{K}$ is not irreducible.
\end{rem}

\begin{pro}\label{pro-xx}
Assume that $K$ is irreducible and
\[
\min_{x\in V} \{K(x,x)\}>0.
\]
Then, for any bijection $g$
of $V$, $\widetilde{K}$ is irreducible and aperiodic, and
$(K_i)_1^\infty$ is merging in relative-sup.
\end{pro}

\begin{pf} By   Example 3.6 of \cite{SZ3}  we have $K_{0,|V|}(x, y)>0$
for all $x,y\in V$. By Corollary  \ref{cor-irr}, this implies that
$\widetilde{K}$ is irreducible aperiodic. By Theorem \ref{th-wave}(1),
we conclude that $(K_i)_1^\infty$ is merging.
\end{pf}

The proof of the proposition above illustrates the surprising
fact that it is not always advantageous to study $\widetilde{K}$
instead of the sequence  $(K_i)_1^\infty$. In  Proposition \ref{pro-xx},
we use the sequence $(K_i)$ to study $\widetilde{K}$!
Indeed, the chain $\widetilde{K}$ seems often difficult to study.
For one thing, $\widetilde{K}$ is not necessarily reversible even if $K$ is.
In general, this means that computing $\widetilde{\pi}$ may be  difficult.
Even when we can compute $\widetilde{\pi}$, it might be difficult
to study the ergodicity of $\widetilde{K}$ from its definition.
Consider, for instance,
the case of cyclic-to-random transposition. In this case,
$\widetilde{\pi}$ is the uniform distribution, but $\widetilde{K}$ is not
invariant under the action of  $S_n$. In other words, the chain driven by
$\widetilde{K}$ is not a random walk on $S_n$. This makes studying
$\widetilde{K}$ and its powers directly rather difficult (and, indeed,
mysterious). The results obtained in \cite{Ga,MPS,SZ} concerning the
cyclic-to-random transposition
chain are essentially obtained by considering the sequence $(K_i)_1^\infty$,
not $\widetilde{K}$ (which, for one thing, does not appear in those papers).

\section{Perturbations of symmetric kernels}\label{sec-pert}

Let $Q$  be a symmetric Markov kernel on a finite set $V$, that is, $Q(x,y)=Q(y,x)$ for
all $x,y\in V$.  This kernel has the uniform distribution $u\equiv 1/|V|$
as its reversible measure.  Fix an $\epsilon\in (0,1)$ and  a set $A\subset V$,
and consider the kernel
\begin{equation}\label{eqn-Ktilde-wave}
K=Q+\Delta_A,
\end{equation}
where $\Delta_A$ is some perturbation kernel such that for
all $x,y\in V$:
\begin{longlist}[(b)]
\item[(a)]$\sum_{z}\Delta_A(x,z)=0$,
\item[(b)]$\Delta_A(x,y)\geq-\epsilon Q(x,y)$  and
\item[(c)]$x\notin A\implies \Delta_A(x,y)=0$.
\end{longlist}

Let $g$ be a permutation of the vertex set $V$ and consider the
sequence $(K_i)_1^\infty$ defined by $K_i(x,y)=K(g^{i-1}x,g^{i-1}y)$.
Set $\widetilde{K}(x,y)=K(x,g^{-1}y)$, as before. Let
$\widetilde{\pi}$ be an invariant probability measure for $\widetilde{K}$
and set
\[
\mu_i(x)=\widetilde{\pi}(g^ix),\qquad  x\in V, i=0,1,2,\ldots.
\]
Define also the symmetric kernel
\[
Q_g(x,y)=Q(g^{-1}x,g^{-1}y).
\]

Consider the following two assumptions on the kernel $\widetilde{K}$:
\begin{longlist}[(A2)]
\item[(A1)](Irreducibility of $\widetilde{K}$) For all $x,y\in V$ there exists an
$n=n(x,y)$
such that $\widetilde{K}^n(x,y)>0$.
\item[(A2)](Aperiodicity of $\widetilde{K}$) There exists a number $N$
such that, for all \vspace{+1pt}
$m\geq N$ and all $x\in V$, $\widetilde{K}^m(x,x)>0$.
\end{longlist}
Recall (see Theorem \ref{th-wave}) that these properties are necessary
for the relative-sup  merging of the sequence
$(K_i)_1^\infty$. In general, it is not obvious at all how they can be checked.
However, if the permutation $g$ is an automorphism of the graph structure on $V$
with edge set $E=\{(x,y) \dvtx K(x,y)>0\}$, then these
properties reduce to the similar properties for $K$ (see Proposition
\ref{pro-kernel-shift}).

The most useful technical result concerning such time inhomogeneous perturbations of
$Q$ is the following comparison lemma.
For more on comparison techniques see \cite{DS-C}.
\begin{lem} \label{lem-Dircomp}
Referring to the above setting, assume that
\begin{equation}\label{MmhypV}
\exists c>0 \qquad  \max_{x\in V}\{\widetilde{\pi}(x)\}\le c \min_{x\in
V}\{\widetilde{\pi}(x)\}.
\end{equation}
Consider the operators $Q_g,\widetilde{K}$ acting respectively on
$\ell^2(u),\ell^2(\widetilde{\pi})$.
Then  the Dirichlet forms $\mathcal E_{Q_g^*Q_g,u}$ of
$Q_g^*Q_g$ on $\ell^2(u)$ and
$\mathcal E_{\widetilde{K}^*\widetilde{K},\widetilde{\pi}}$ of
$\widetilde{K}^*\widetilde{K}$ on $\ell^2(\widetilde{\pi})$ satisfy
\begin{equation}\label{Dircomp}
\mathcal E_{Q_g^*Q_g,u}(f,f)\le  \frac{c}{(1-\epsilon)^2}
\mathcal E_{\widetilde{K}^*\widetilde{K},\widetilde{\pi}}(f,f)
\end{equation}
for any function $f$ defined on $V$.
\end{lem}

\begin{pf}Working on $\ell^2(\widetilde{\pi})$ and $\ell^2(u)$, respectively,
we compare the kernel
$\widetilde{K}^*\widetilde{K}$ to
the kernel $Q_g^*Q_g$, that is,  $Q^*Q$ moved by $g^{-1}$.  Write
\begin{eqnarray*}
\widetilde{\pi}(x)
\widetilde{K}^*\widetilde{K}(x,y)
&\ge & \frac{1}{c}\sum_z u(z)K(z,g^{-1}x)K(z,g^{-1}y)\\
&\ge & \frac{(1-\epsilon)^2}{c}\sum_zu(z)Q(z,g^{-1}x)Q(z,g^{-1}y)\\
&=& \frac{(1-\epsilon)^2}{c} u(x)Q_g^* Q_g( x,y).
\end{eqnarray*}
The third line uses the fact that for any $z$, $u(g^{-1}z)= u(z)= 1/|V|$.
\end{pf}

The importance of this lemma comes from the fact that $Q_g$
is simply $Q$ transported by $g^{-1}$ and thus has the same properties as
$Q$. For instance, $Q_g$ has the same eigenvalues and singular values as $Q$
(the eigenvectors of $Q_g$ are the eigenvectors
of $Q$ transported by $g^{-1}$, etc.).  Similarly, $Q_g$
satisfies the same Nash and logarithmic Sobolev inequalities on $\ell^2(u)$
as $Q$ itself. By Lemma \ref{lem-Dircomp}, these properties will be
transferred to
$(\widetilde{K},\widetilde{\pi})$. The following two propositions and assorted
remarks are based on this observation.

\begin{pro}\label{cor-wave}
Referring to the above setting, assume that (\ref{MmhypV}) holds, that is,
\[\max_{x\in V}\{\widetilde{\pi}(x)\}\le c \min_{x\in V}\{\widetilde{\pi}(x)\}.
\]
Let $\sigma_1$ be the second largest singular value of $Q$ on $\ell^2(u)$.
Then the second largest singular value $\widetilde{\sigma}_1$ of $\widetilde{K}$ on
$\ell^2(\widetilde{\pi})$ is bounded by
\[
\widetilde{\sigma}_1\le 1-\frac{(1-\epsilon)^2}{c^2}(1-\sigma_1).
\]
Furthermore by Theorem \ref{th-wave} we obtain
\[\max_{x,z\in V}\biggl \{ \bigg|\frac{K_{0,n}(x,z)}{\mu_n(z)}-1 \bigg| \biggr\}\le c|V|
 \biggl(1-\frac{(1-\epsilon)^2}{c^2}(1-\sigma_1) \biggr)^n.
 \]
\end{pro}

\begin{rem}If instead of using $\sigma_1$ we use the logarithmic Sobolev constant
$l(Q^*Q)$ of $Q^*Q$ (see \cite{DS-L,SZ4} for the definition;
we follow the notation of \cite{SZ4});
then we get
\[
l(\widetilde{K}^*\widetilde{K})\ge \frac{(1-\epsilon)^2}{c^2} l(Q^*Q).
\]
In cases where a good estimate on $l(Q^*Q)$ is known, this can, potentially,
improved upon the merging bound stated in the corollary above. See
\cite{DS-L,SZ4}.
\end{rem}

In the next corollary, we make use of one of the main results of
\cite{DS-N,SZ4}
which concerns the use of the Nash inequalities.
In applications, the constants $c$, $c_1$, $C_1$, $D$ appearing in the statement below are indeed taking fixed values
whereas the parameter $T$ grows with the size of the underlying state space.
It is, in general, equal to the square of the diameter of the state space
$V$ equipped with the graph structure induced by the symmetric kernel $Q$.
For an introduction to the use of Nash inequality in the study of ergodic
Markov chains, see \cite{DS-N}.

\begin{pro}\label{cor-wave-Nash}
Referring to the above setting, assume that there are constants
$c,c_1,C_1,D\in (0,\infty)$ and a parameter $T>1$ such that:
\begin{itemize}
\item  Condition (\ref{MmhypV}) holds, that is,
\[\max_{x\in V}\{\widetilde{\pi}(x)\}\le c \min_{x\in V}\{\widetilde{\pi}(x)\}.
\]
\item The second largest singular value $\sigma_1(Q)$ of $Q$ on $\ell^2(u)$ satisfies
\[
\sigma_1(Q)\le 1- \frac{c_1}{ T}.
\]
\item  The kernel $Q$ satisfies the Nash inequality (all norms are w.r.t. $u$)
\[
\forall  f \dvtx V\rightarrow V \qquad  \|f\|_2^{2 +1/D}
\le C_1T \biggl(\mathcal E_{Q^*Q}(f,f)+\frac{1}{T}\|f\|_2^2 \biggr)\|f\|_1^{1/D}.
\]
\end{itemize}
Then, for any $n> 2T$ and $x,z\in V$, we have
\[ \bigg|\frac{K_{0,n}(x,z)}{\mu_n(z)}-1 \bigg|\le
  \biggl(\frac{16(1+4D)C_1c^{2+3/(2D)}}
{(1-\epsilon)^2} \biggr)^{2D} e^{- 2c_1(1-\epsilon)^2(n-2T)/c^2T}.
\]
\end{pro}

\begin{pf}
Let $u\equiv 1/|V|$. For any function $f \dvtx V\ra V$ we have
$\mathcal E_{Q^*Q,u}(f,f)=
\mathcal E_{Q_{g}^*Q_g,u}(f\circ g^{-1},f\circ g^{-1})$
and $\|f\|_p=\|f\circ g^{-1}\|_p$ for $p=1,2$.\vspace*{1pt} Thus
$(\mathcal E_{Q^*_gQ_g},u)$ satisfies the same Nash inequality
as $(\mathcal E_{Q^*Q},u)$. By Lemma \ref{lem-Dircomp} and (\ref{MmhypV}),
this yields the
Nash inequality,
\[
\|f\|_{\ell^2(\widetilde{\pi})}^{2 +1/D}
\le \frac{C_1Tc^{2+3/(2D)}}{(1-\epsilon)^2}
\biggl (\mathcal E_{\widetilde{K}^*\widetilde{K},\widetilde{\pi}}
(f,f)+\frac{1}{T}\|f\|_{\ell^2(\tilde{\pi})}^2 \biggr)\|f\|_{\ell^1(\widetilde{\pi})}^{1/D}
\]
for $(\mathcal E_{\widetilde{K}^*\widetilde{K}},\widetilde{\pi})$.
The desired result now follows by applying Propositions \ref{pro-kernel-shift},
\ref{cor-wave} and  the results of \cite{DS-N}. (See also
Theorem 2.5 of \cite{SZ4}.)
\end{pf}

Observe that the conclusion can be rephrased by saying that, under the hypotheses made,
the time inhomogeneous chain driven by $(K_i)_1^\infty$ has a relative-sup merging time
at most of order $T$. This will be illustrated below in concrete examples.

Assuming (as is natural) that we understand well the finite Markov chain driven
by the symmetric kernel $Q$, the main difficulty that remains in studying the
time inhomogeneous chain $(K_i)_1^\infty$ considered in this section is to verify
the condition~(\ref{MmhypV}) for some (explicit)  constant $c$.
The following lemma is useful in this regard.

\begin{lem}\label{lem-wave-minmax} Assume that $\widetilde{\pi}\neq u$ and that
 $\widetilde{K}$  satisfies the irreducibility
condition \textup{(A1)} above.  Let $M=\max_x\{\widetilde{\pi}(x)\}$ and
$m=\min_x\{\widetilde{\pi}(x)\}$. Let
\[
A^*_+=\biggl\{x\in V \dvtx \sum_y \widetilde{K}(y,x)>1\biggr\}, \qquad
 A^*_-=\biggl\{x\in V \dvtx \sum_y \widetilde{K}(y,x)<1\biggr\}.
\]
 Then there are points
$x_+ \in A^*_+,x_-\in A^*_-$ such that $\widetilde{\pi}(x_+)=M$,
$\widetilde{\pi}(x_-)=m$.
\end{lem}

\begin{pf} Let
$B=\{z \dvtx \sum_y \widetilde{K}(y,z)=1\}.$
Let $x\in V$ be a point such that
$\widetilde{\pi}(x)=M$.
Then we must have $\sum_y\widetilde{K}(y,x)\ge 1$. If $\sum_y\widetilde{K}(y,x)>1$,
we are done. Otherwise, $x\in B$ and we must have $\widetilde{\pi}(y)=M$ for all $y$
such that $\widetilde{K}(y,x)>0$. Either one of these points $y$ satisfies
$\sum_z\widetilde{K}(z,y)>1$ and we are done, or we repeat the argument. Since $\widetilde{K}$
satisfies (A1) and $\widetilde{\pi}\neq u$, this process necessarily yields a point~$x_+$
such that $\widetilde{\pi}(x)=M$ and $x_+\notin B$. Of course,
we must then have $x_+\in A^*_+$. The same line of reasoning proves the existence of the
desired point $x_-\in A^*_-$.
\end{pf}

\begin{rem}\label{rem-A*}
Note that $A^*_+, A^*_-$ are contained in the ``$\widetilde{K}$-boundary'' of $A$, that is in the
set $A^*=\{z: \exists  y\in A,  \widetilde{K}(y,z)>0\}$. Indeed, if $x\notin A^*$ then
\[
\sum_y \widetilde{K}(y,x)= \sum_y Q(y,g^{-1}x)=\sum_y Q(g^{-1}x,y)=1.
\]
{\smallskipamount=0pt
\vspace*{-6pt}
\begin{longlist}[(b)]
\item[(a)] If we can find $n_0$ such that
$\inf\{\widetilde{K}^{n_0}(x,y) \dvtx x,y\in A^*\}>\delta>0$,
then since $\widetilde{\pi}= \widetilde{\pi}\widetilde{K}^{n_0}$, one obtains
$\widetilde{\pi}(x_+)=\max\{\widetilde{\pi}\}\le \delta^{-1} \min\{\widetilde{\pi}\}   =\delta^{-1}\widetilde{\pi}(x_-).$
Unfortunately, the nature of the kernel $\widetilde{K}$ makes it difficult to find a suitable $n_0$.
\item[(b)] A variation on this idea is as follows. Assume that, for any $(x,y)\in A^*_+\times A^*_-$, we
can find an element $b=b(x,y)$ such that
\[
\frac{\widetilde{K}(b,x)}{1-\sum_{z\neq b} \widetilde{K}(z,x)}\in (0,\infty) \quad \mbox{and}\quad
\frac{1-\sum_{z\neq b}\widetilde{K}(z,y)}{\widetilde{K}(b,y)}\in (0,\infty).
\]
Then for $x,y\in A^*_+\times A^*_-$ such that $\widetilde{\pi}(x)=M$ and $\widetilde{\pi}(y)=m$ as
defined in Lemma~\ref{lem-wave-minmax} we have
\[
\widetilde{\pi}(x)\le  \biggl(\frac{\widetilde{K}(b,x)(1-\sum_{z\neq b}\widetilde{K}(z,y))}
{\widetilde{K}(b,y) (1-\sum_{z\neq b} \widetilde{K}(z,x))}
 \biggr)\widetilde{\pi}(y).
 \]
 This gives
$\max\{\widetilde{\pi}\}\le C\min\{\widetilde{\pi}\}$ with
\[
C=\max_{(x,y)\in A^*_+\times A^*_-}
\biggl \{\frac{\widetilde{K}(b,x)(1-\sum_{z\neq b}\widetilde{K}(z,y))}
{\widetilde{K}(b,y) (1-\sum_{z\neq b} \widetilde{K}(z,x))} \biggr\}.
\]
\end{longlist}}
\end{rem}

Note that $C$ depends on the choice of the $b(x,y)$ for each $(x,y)\in A_+^*\times A_-^*$.
Different choices
of allowed $b$s may yield a different constant $C$.  If the location of $\max \widetilde{\pi}$
and $\min \widetilde{\pi}$
can be determined, then there is no need to calculate $C$ over all $A_+^*\times A_-^*$.
Examples using this remark are in the next two sections.

\section{Cyclic edge perturbation on the circle} \label{ex-twc}
This section examines some examples of a moving wave on the circle graph.
On the circle graph on $N=2 l+1$ vertices and for $\epsilon>0$ fixed,
let $K$ be the reversible Markov kernel corresponding to putting weight $1$
on all edges except the $(0,1)$ edge which has weight $1+\epsilon$. Hence
\begin{equation}\label{def-Kcircle}
K(x,y)=
\cases{
0, & \quad if  $|x-y|\neq 1$,\cr
1/2, & \quad if $|x-y|=1 \mbox{ and } x \notin \{0,1\}$,\cr
(1+\epsilon)/(2+\epsilon), &\quad if $(x,y)\in \{(0,1),(1,0)\}$,\cr
1/(2+\epsilon),& \quad if $(x,y)\in \{(0,-1),(1,2)\}$.
}
 \end{equation}
 This has reversible measure
\[
\pi(x)= \cases{
1/(N+\epsilon),&\quad if $x\neq 0,1$,\cr
(1+\epsilon/2)/(N+\epsilon), &\quad if $x=0,1$.
}
\]
Note that this can be written as a perturbation (see Section \ref{sec-pert})
of the symmetric kernel $Q$ of simple random walk, $Q(x,y)=1/2$ if $|x-y|=1$
and $Q(x,y)=0$ otherwise. The perturbation set
$A$ is $A=\{0,1\}$  and $\Delta_A=0$ except for the following values:
\[
\Delta_A(0,1)=\Delta_A(1,0)=\epsilon/(4+2\epsilon),\qquad
\Delta_A(0,-1)=\Delta_A(1,2)=-\epsilon/(4+2\epsilon).
\]
Because $N=2l+1$ is odd, the chain driven by $Q$ is ergodic with
relative-sup mixing time of order $N^2$. Its singular values (i.e.,  eigenvalues) on $\ell^2(u)$ are
\[
\cos \biggl(\frac{2\pi j}{N} \biggr),\qquad   j=0, 1,\ldots, N-1.
\]
In particular, the second largest is attained at $j=(N-1)/2$ and equals
\begin{equation}\label{circ1}
\beta_1=\cos\frac{\pi}{N}.
\end{equation}
Moreover, $Q$ satisfies the Nash inequality
\begin{equation}
\label{circ2} \qquad\forall f \dvtx V\rightarrow V \qquad   \|f\|_2^6\le 2^7N^2
\biggl(\mathcal E_{Q^*Q}(f,f)+\frac{1}{4(N+1)^2}\|f\|_2^2 \biggr)\|f\|_1^{4}   .
\end{equation}
See, for example, Theorem 5.2 and Lemma 5.3 in \cite{DS-N}.

We will investigate the general construction described  earlier
based on the kernel $K$ above and  various  bijections including
$x\mapsto x-1$ and $x\mapsto x+2$. In these two cases, we prove a merging time
estimate of the type
\[
T_\infty(\eta)\le C(\epsilon)N^2(1+\log_+1/\eta) \qquad \forall \eta>0
\]
for the associated periodic time inhomogeneous chain, but
there are interesting differences in the analysis of the two chains.

\begin{figure}

\includegraphics{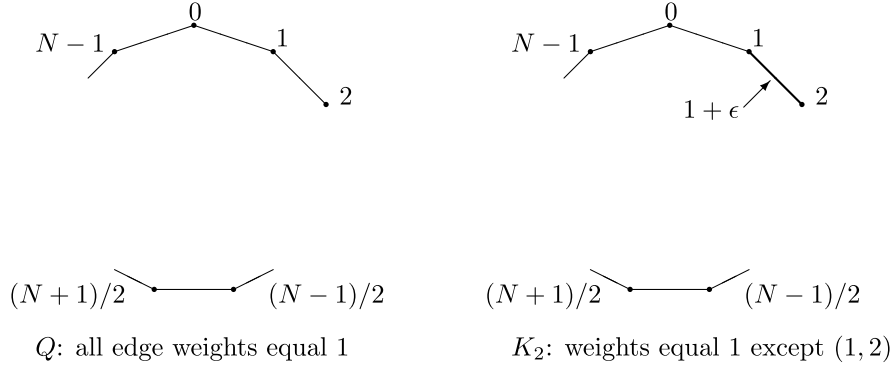}

\caption{The cycling edge perturbation of $Q$.}\label{fig3}
\end{figure}

First, consider $g(x)=x-1$. Then $K_i$ is the reversible kernel
corresponding to
putting weight $1+\epsilon$ on the edge $(i-1,i) \mod N$.
The graphs\vadjust{\goodbreak} for $Q$ and $K_2$ are given in Figure~\ref{fig3}.
The kernel
$\widetilde{K}(x,y)=K(x,g^{-1}y)$ is given by
\[
\widetilde{K}(x,y)=
\cases{ 0, & \quad if $y\notin\{x,x-2\}$,\cr
1/2, &\quad if  $y\in \{x,x-2\}$   and   $x\notin \{0,1\}$,\cr
(1+\epsilon)/(2+\epsilon),&\quad if $(x,y)\in \{(0,0),(1,-1)\}$,\cr
1/(2+\epsilon),  &\quad if $(x,y)\in \{(0,-2),(1,1)\}$.
}
\]
A simple calculation shows that
$\widetilde{\pi}$ is constant away from $0,1$ and  that
\[
\widetilde{\pi}(x)=
\cases{
{2(1+\epsilon)}/{(\epsilon^2+2N\epsilon+2N)}, &\quad if $x\neq
0,1$,\cr
 {(\epsilon+1)(\epsilon+2)}/{(\epsilon^2+2N\epsilon+2N)}, &\quad if
 $x=0$,\cr
{(\epsilon+2)}/{(\epsilon^2+2N\epsilon+2N)}, &\quad if $x=1$.
}
\]
This proves $c$-stability of the sequence $(K_i)_1^\infty$ with respect to
$\mu_0=\widetilde{\pi}$ with $c=1+\epsilon$. This distribution yields the
wave $\mu_i(x)=\widetilde{\pi}(g^i x)$ created  by the time
inhomogeneous Markov chain driven by $(K_i)_1^\infty$.

Using Proposition \ref{cor-wave} and (\ref{circ1}),
this proves that the  relative-sup merging time for the sequence $(K_i)_1^\infty$ is
bounded by $T_\infty(\eta)\le C(\epsilon)N^2(\log N+\log_+ 1/\eta)$.
An improved result showing relative-sup merging in time of order $N^2$ is obtained
using Proposition \ref{cor-wave-Nash} and the Nash inequality (\ref{circ2}) of the circle graph.

\begin{figure}
    
\includegraphics{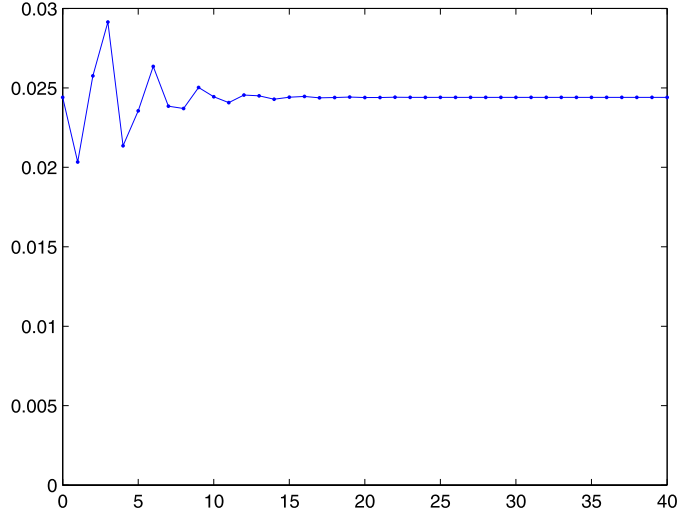}

  \caption{$\tilde{\pi}$ for $N=41$ and $\varepsilon=1$.}\label{fig4}
\end{figure}

Let us now consider what happens if we choose $g(x)= x+2$.
In terms of the sequence $K_i$, this means that $K_i$ now has the same
perturbation as $K$ but at the edge $(-2i,-2i+1) \mod N$.
The kernel $\widetilde{K}$ is given by
\[
\widetilde{K}(x,y)= \cases{
 0,&\quad if
$y-x\notin \{1,3\}$,\cr
1/2, & \quad if $y-x\in\{1,3\}$   and  $x\notin \{ 0,1\}$,\cr
(1+\epsilon)/(2+\epsilon),& \quad if $(x,y)\in\{(0,3),(1,2)\}$,\cr
1/(2+\epsilon), & \quad if $(x,y)\in \{(0,1),(1,4)\}$.
}
\]
Contrary to what happens with $g:x\mapsto x -1$, in the present case, there
is no simple formula for $\widetilde{\pi}$ (in particular,
$\widetilde{\pi}$ is not  constant away from the perturbation).
Figure~\ref{fig4} presents a simulation of the stationary measure $\widetilde{\pi}$
for $N=41$ and $\epsilon=1$.

However, it is easy to see from the linear equations defining
$\widetilde{\pi}$ (i.e., from Lemma~\ref{lem-wave-minmax})
that $\max\{\widetilde{\pi}\}$ must be attained at either
$2$ or $3$,
and $\min \{\widetilde{\pi}\}$ must be attained at either $1$ or $4$.
Suppose they are attained at $2$ and $1$. As
\[\widetilde{\pi}(2)=  \biggl(\frac{1+\epsilon}{2+\epsilon} \biggr)\widetilde{\pi}(1) +
\frac{1}{2}\widetilde{\pi}(4)
\]
we must have
\[
\widetilde{\pi}(2)\le  \biggl(\frac{1+\epsilon}{1+\epsilon/2} \biggr) \widetilde{\pi}(1).
\]
Suppose instead the max and min  are attained at  $2$ and $4$. Then, the same
equation gives
\[
 \biggl(1-\frac{1+\epsilon}{2+\epsilon} \biggr) \widetilde{\pi}(2)\le
\frac{1}{2}\widetilde{\pi}(4),
\]
that is,
\[
\widetilde{\pi}(4)\ge \biggl( \frac{1}{1+\epsilon/2}  \biggr)\widetilde{\pi}(2).
\]
The case where the max and min are attained at $3$ and $2$ is treated similarly. The
remaining case where the max and min are attained at $3$ and $1$ is slightly different
because there is no direct relation between
$\widetilde{\pi}(3)$ and $\widetilde{\pi}(1)$. However, the same line of reasoning
yields
\[
\widetilde{\pi}(3)\le  \biggl(\frac{1+\epsilon}{1+\epsilon/2}  \biggr)\widetilde{\pi}(0)
\quad\mbox{and}\quad \widetilde{\pi}(0)\le (1+\epsilon/2) \widetilde{\pi}(1).
\]
 This shows that
\begin{equation}\label{shift2}
\max\{\widetilde{\pi}\}\le (1+\epsilon)
\min\{\widetilde{\pi}\}.
\end{equation}
Because of this and Corollary \ref{mv-stab},
the sequence $(K_i)_1^\infty$ is $(1+\epsilon)$-stable
with respect to $\widetilde{\pi}$.  Applying Proposition \ref{cor-wave-Nash}
and (\ref{circ2}) yield again a relative merging time  of order $N^2$ for the
sequence $(K_i)_1^\infty$. The following theorem records this result in more general form.

\begin{theo}\label{thm-RS-circ}
Let $V_N=\{0,\ldots,N\}$. Fix $\epsilon>0$  and let $K$ be as
in (\ref{def-Kcircle}).  Fix a permutation $g=g_N$ of $V_N$ and
let $K_i$, $\widetilde{K},\widetilde{\pi}, \mu_i$
be associated to $K,g$ as  in \textit{Section~\ref{sec-wave}}.
Assume that  there exists $c\ge 1$ such that
\begin{equation}\label{Mmhyp}
\max_{x\in V_N}\{\widetilde{\pi}(x)\}\le c \min_{x\in
V_N}\{\widetilde{\pi}(x)\}.\end{equation}
Then there is a constant $C(\epsilon,c)$ such that the relative-sup merging time
for $(K_i)_1^\infty$ is bounded by
\[
T_\infty(\eta)\le C(\epsilon,c) N^2(1+\log_+1/\eta).
\]
\end{theo}

\begin{rem}\label{rem-circ-map}
For which permutations $g$ of the set $V_N=\{0,\ldots, N\}$
does  the conclusion of the theorem above hold?  According to the theorem,
it suffices to check that condition (\ref{Mmhyp}) is satisfied.
For instance, (\ref{Mmhyp}) is satisfied if $g(x)=x-1$ or $g(x)= x+ 2$
[in fact, by symmetry,
for $g(x)=x\pm 1$, $g(x)=x\pm 2$].
It is very plausible that (\ref{Mmhyp}) is always satisfied, whatever the
permutation $g$ is.
However, this does not follow directly from an argument similar to the one used for
$g(x)=x-1$ and $g(x)=x+2$. In fact, the argument already fails
miserably for $g(x)=x+3$. The reader may want to convince herself of that.
In general, we want to compare the min and max of $\widetilde{\pi}$. It is easy to
see that
the max is attained at either $g(0)$ or $g(1)$ and the min at
either $g(-1)$ or $g(2)$.  The case where the max and min are attained
at either $(g(0),g(2))$ or $(g(1),g(-1))$ can be treated as above because
the values of $\widetilde{\pi}$ at $g(0),g(2)$
[resp., at $g(-1),g(1)$] are both related to the value at $1$ (resp., $0$).
But, in  the other cases, it becomes much more tricky to compare the max and min
without further hypotheses.
\end{rem}

Let $P$ be the lazy version of the kernel defined in (\ref{def-Kcircle}) with
\begin{equation}\label{def-Qcircle}
P(x,y)= \cases{
1/2, &\quad if $x=y$,\cr
1/4, &\quad if $|x-y| = 1$ and $x\neq\{0,1\}$,\cr
(1+\epsilon)/2(2+\epsilon), &\quad if $(x,y)\in\{(0,1),(1,0)\}$,\cr
1/2(2+\epsilon), &\quad if $(x,y)\in\{(0,-1),(1,2)\}$,\cr
0, & \quad otherwise.}
\end{equation}
Let $g$ be any permutation of the set $V_N=\{0,\ldots, N\}$, and define
$P_i(x,y)=P(g^{i-1}x,g^{i-1}y)$  for all $i=1,2,\ldots$  and
$\widetilde{P}(x,y)=P(x,g^{-1}y)$.
In this case, we can show that condition (\ref{MmhypV}) holds which implies a
relative-sup merging time of order $N^2$ for any permutation $g$.

\begin{theo}
Let $V_N=\{0,\ldots,N\}$. Fix $\epsilon>0$  and let $P$ be as
in (\ref{def-Qcircle}).  Fix a permutation $g=g_N$ of $V_N$ and
let $P_i$, $\widetilde{P},\widetilde{\pi}, \mu_i$
be associated to $P,g$ as  in \textit{Section~\ref{sec-wave}} (replacing $K$ by $P$).
Then
\begin{equation}\label{max/min}
\max_{x\in V_N}\{\widetilde{\pi}(x)\}\le (1+\epsilon)
\min_{x\in V_N}\{\widetilde{\pi}(x)\}.\end{equation}
Furthermore,
there is a constant $C(\epsilon)$ such that the relative-sup merging time
for $(P_i)_1^\infty$ is bounded by
\[
T_\infty(\eta)\le C(\epsilon) N^2(1+\log_+1/\eta).
\]
\end{theo}

\begin{pf}  By Proposition \ref{cor-wave-Nash} and (\ref{circ1})--(\ref{circ2}),
it suffices to prove (\ref{max/min}).
Fix a permutation $g=g_N$ of $V_N=\{0,\ldots,N\}$. The kernel $\widetilde{P}$ is given by
\[
\widetilde{P}(x,y)= \cases{
1/2, &\quad if  $x = g^{-1}y$,\cr
1/4, &\quad if $|x-g^{-1}y| = 1$ and $x\neq\{0,1\}$,\cr
(1+\epsilon)/2(2+\epsilon), &\quad if
$(x,g^{-1}y)\in\{(0,1),(1,0)\}$,\cr
1/2(2+\epsilon), &\quad if $(x,g^{-1}y)\in\{(0,-1),(1,2)\}$,\cr
0, &\quad  otherwise.
}
\]
By Lemma \ref{lem-wave-minmax}, the maximum value of $\widetilde{\pi}$
is attained at either $g(0)$ or $g(1)$ and the minimum at $g(-1)$ or $g(2)$. Moreover,
\begin{eqnarray*}
\widetilde{\pi}(g(-1))&=&\frac{\widetilde{\pi}(-1)}{2}+\frac{\widetilde{\pi}(-2)}{4}
+\frac{\widetilde{\pi}(0)}{2(2+\epsilon)},\\
  \widetilde{\pi}(g(2))&=&\frac{\widetilde{\pi}(2)}{2}+\frac{\widetilde{\pi}(3)}{4}
+\frac{\widetilde{\pi}(1)}{2(2+\epsilon)},
\\
\widetilde{\pi}(g(0))&=&\frac{\widetilde{\pi}(0)}{2}+\frac{\widetilde{\pi}(-1)}{4}
+\frac{(1+\epsilon)\widetilde{\pi}(1)}{2(2+\epsilon)},\\
\widetilde{\pi}(g(1))&=&\frac{\widetilde{\pi}(1)}{2}+\frac{\widetilde{\pi}(2)}{4}
+\frac{(1+\epsilon)\widetilde{\pi}(0)}{2(2+\epsilon)}.
\end{eqnarray*}
Note that for any of the four possible max/min pairs, the max and min values
can be both compared via the equations  above to either  $\widetilde{\pi}(0)$
or $\widetilde{\pi}(1)$. See Remark~\ref{rem-A*}(b).
For instance,  suppose the max/min pair is $(g(0),g(-1))$. Then
\[
\widetilde{\pi}(g(0))\leq \frac{4+2\epsilon}{4+\epsilon}\widetilde{\pi}(0)
 \quad\mbox{and}\quad
\widetilde{\pi}(0)\leq\frac{2+\epsilon}{2}\widetilde{\pi}(g(-1)).
\]
Hence,
\[
\widetilde{\pi}(g(0))\leq\frac{(2+\epsilon)^2}{4+\epsilon}\widetilde{\pi}(g(-1)).
\]
The other cases are similar, and
it follows that $\max\{\widetilde{\pi}\}\leq (1+\epsilon)\min\{\widetilde{\pi}\}$.
\end{pf}

\section{Further examples: Single point perturbations}

In the next two examples, we consider perturbations of a symmetric kernel as
described in Section \ref{sec-pert} but with $A=\{o\}$ for some $o\in V$, that is,
the perturbation occurs at a single point.
In the second example, we make an additional assumption on the structure of the
perturbation. In these cases, we are able to obtain easily applicable bounds.
\begin{example}\label{ex-wave-vertex} Let $Q$ be be a symmetric kernel as in
Section \ref{sec-pert}. Fix $\epsilon\in(0,1)$, and let
$K=Q+\Delta_o$ where $\Delta_o=\Delta_{\{o\}}$ satisfies
\[
-\epsilon Q(o,y)\leq \Delta_o(o,y),\qquad  \sum_y\Delta_o(o,y)=0  \quad\mbox{and}\quad
\Delta_o(x,y)=0 \qquad\mbox{if } x\neq o.
\]
Note that  $K(x,y)\geq (1-\epsilon)Q(x,y)$, and
$K$ satisfies the properties (a)--(c) listed at the beginning of Section \ref{sec-pert}.  Fix a permutation $g$ of $V$ and assume that
$\widetilde{K}$ is irreducible.
Then Lemma \ref{lem-wave-minmax} says that the min and max of $\widetilde{\pi}$ are
attained respectively  on $A^*_+,A^*_-$ and Remark \ref{rem-A*}(b) gives
\begin{equation}\label{stabtheta}
\max_{x\in V}\{\widetilde{\pi}\}\le C
\min_{x\in V}\{\widetilde{\pi}\},
\end{equation}
where
\begin{eqnarray*}
C&=&
\max_{(x,y)\in A^*_+\times A^*_-}\biggl \{\frac{\widetilde{K}(o,x)(1-\sum_{z\neq o}
\widetilde{K}(z,y))}{\widetilde{K}(o,y) (1-\sum_{z\neq o} \widetilde{K}(z,x))} \biggr\}\\
&\le &
\max_{x\in A^*_+} \biggl\{\frac{K(o,g^{-1}x)}{ (1-\epsilon)Q(o,g^{-1}x)} \biggr\}\\
&= &
\frac{1}{(1-\epsilon)\theta}, \qquad  \theta= \max_{x\in A^*_+} \biggl\{\frac{K(o,g^{-1}x)}{Q(o,g^{-1}x)} \biggr\}.
\end{eqnarray*}

Equation (\ref{stabtheta}) and Proposition \ref{cor-wave} now imply that
the relative-sup $\eta$ merging time of the sequence $(K_i)_1^{\infty}$ is  at most
\begin{equation} \label{mergetheta}
\frac{D}{1-\sigma_1} (\log |V| +\log_+1/\eta ),
\end{equation}
where $\sigma_1$ is the second largest singular value of the kernel $Q$ on
$\ell^2(u)$, and $D=D(\epsilon,\theta)$
is a constant  that depends only on $\epsilon\in (0,1)$ and
$\theta$ (the constant $D$ can easily be made explicit).
\end{example}

\begin{example}[(Perturbation of expander graphs)] \label{exp1}
Fix an integer $r$ and
consider a sequence $\mathcal G_N=(V_N,E_N)$ of regular graphs with vertex set
$V_N$ of size $|V_N|$ tending to infinity
and symmetric edge set $E_N\subset V_N\times V_N$
with $(x,x)\in E_N$ for all $x\in V_N$.
On each graph, consider the symmetric  Markov kernel $Q=Q_N$ corresponding
to the simple random walk on $\mathcal G_N$.
Hence, $Q_N(x,y)=1/r$
if $(x,y)\in E_N$ and $Q_N(x,y)=0$ otherwise.
Let $\sigma_1(N)$ be the
second largest singular value of $Q_N$ on $\ell^2(u_N)$ where $u_N$ is the
uniform probability measure on $V_N$. Assume
that there is a constant $a\in (0,1)$ such that
\begin{equation}\label{expdr}
\forall N \qquad  1-\sigma_1(N)\ge a.
\end{equation}
This property is a strong form of the property  that defines the so-called
expander graphs (see, e.g., \cite{Lub1,Lub2} and the references therein).

Fix an origin $o=o_N$ in $V_N$ and consider a perturbation $K_N$ of
$Q_N$ as  in Example \ref{ex-wave-vertex}.
Fix also a bijection $g_N \dvtx V_N\ra V_N$. For each $N$,
consider the time inhomogeneous chain on $V_N$ driven by
$(K_{N,i})_1^\infty$ where $K_{N,i}(x,y)=K_N(g_N^{i-1}x,g_N^{i-1}y)$.
In this situation, (\ref{mergetheta}) yields  merging for the sequence
$(K_{N,i})_1^\infty$ in order $\log |V_N|$ steps, uniformly in $N$.
Note that this result requires the degree $r$ of the graph to be fixed
(or, at least, bounded from above, uniformly in $N$).
\end{example}

\begin{example}\label{ex-wave-vertex'}
Here we strengthened the hypotheses and the conclusion in the previous example.
Namely, we assume that there exists $\delta\in (0,1-Q(o,o))$ such that
\begin{equation}\label{vertex'}
\qquad 0<\Delta_o(o,o)\leq\delta,\qquad
-\delta \biggl(\frac{Q(o,y)}{1-Q(o,o)} \biggr)\leq \Delta_o(o,y)< 0\qquad
\mbox{if }y\neq o,
\end{equation}
and
\[
\Delta(x,y)=0 \qquad\mbox{if }x\neq o.
\]
Set
\begin{equation}\label{eqn-epsilon-wave}
\epsilon= \frac{\delta}{1-Q(o,o)}.
\end{equation}
A careful analysis of this example yields a much improved estimate
for $c$-stability and the relative sup merging time when compared to the previous
example. The difference lies in the fact that the perturbation is positive only at
$o$.
\begin{lem} \label{lem-mMstab}
Assume that $\widetilde{K}$ is irreducible.
Let $m=\min_x\{\widetilde{\pi}(x)\}$ and $M=\max_x\{\widetilde{\pi}(x)\}$. We have that
$\widetilde{\pi}(o)=M$ and for $\epsilon$ as in (\ref{eqn-epsilon-wave})
\[
m\geq (1-\epsilon) \widetilde{\pi}(o).
\]
\end{lem}

\begin{pf}
Lemma \ref{lem-wave-minmax} tells us that $M=\widetilde{\pi}(o)$ and that
there exists $m=\widetilde{\pi}(x_-)$ for some $x_-$ with $\widetilde{K}(o,x_-)>0$.
Further,
\begin{eqnarray*}
\widetilde{\pi}(x_-)&=&\sum_x\widetilde{\pi}(x)\widetilde{K}(x,x_-)\\
&\geq&\widetilde{\pi}(o)\widetilde{K}(o,x_-)+\widetilde{\pi}(x_-)
\sum_{x\neq o}Q(x,g^{-1}x_-)\\
&\geq& (1-\epsilon)\widetilde{\pi}(o)Q(o,g^{-1}x_-)+\widetilde{\pi}(x_-)\bigl(1-Q(o,g^{-1}x_-)\bigr).
\end{eqnarray*}
So we get $\widetilde{\pi}(x_-)\geq (1-\epsilon)\widetilde{\pi}(o)$ as desired.
\end{pf}
\end{example}

\begin{example} Let $\mathcal G_N=(V_N,E_N)$ be
a sequence of regular expander graphs as in  Example \ref{exp1} but
with degree $r_N\ge 3$ that might depend on $N$.
Fix $\delta \in (0,2/3)$ and bijections $g_N \dvtx V_N\ra V_N$.
Consider a perturbation $K_N$ of the simple random walk $Q_N$ on $\mathcal G_N$
as in Example \ref{ex-wave-vertex'}. The constant
$\epsilon$ at (\ref{eqn-epsilon-wave}) is
$\epsilon_N= \delta (r_N/(r_N-1))< 3\delta/2$ and the measure $\widetilde{\pi}_N$
satisfies
\[\max_{V_N}\{\widetilde{\pi}_N\}\le (1-3\delta/2)^{-1} \min_{V_N}
\{\widetilde{\pi}_N\}.
\]  It follows from this and Proposition \ref{cor-wave}
that the associated sequence of perturbed kernels
$(K_{N,i})_1^\infty$ merges in order $\log |V_N|$ steps.
\end{example}

\begin{example}[(Sticky permutation)]
The following is a particular case of  Example \ref{ex-wave-vertex'}.
It is treated in more detail in \cite{SZ4}. On $V=S_n$, the symmetric group,
let
\[
Q(x,y)=
 \cases{
1/2n, &  \quad if   $y=x(1,j)$,   $j\in \{2,\ldots,n\}$,\cr
(n+1)/(2n), & \quad if $x=y$.\cr
0, &\quad   otherwise.
}
\]
This is the kernel of the lazy version of  the random walk called
``transpose top and random.'' Fix a permutation $\rho_n\in S_n$,
$\delta\in (0,(n-1)/(2n))$
and let
\[
K(x,y)= \cases{
Q(x,y),& \quad if $x\neq\rho_n$,\cr
Q(x,y)+\delta, & \quad if $x=y=\rho_n$,\cr
Q(x,y)-\delta/(n-1), & \quad if $x=\rho_n$ and $y=x(1,j)$\cr
&\qquad for  $j\in \{2,\ldots,
n\}$.
}
\]
In words, $K$ is obtained from $Q$ by  adding extra holding
probability at $\rho_n$,  making~$\rho_n$ ``sticky.''
Next, if $\sigma$ is the cycle $(1,\ldots,n)$, let
\[
K_i(x,y)= K(\sigma^{i-1} x \sigma^{-i+1},\sigma^{i-1} y
\sigma^{-i+1}).
\]
Hence $K_i$ is $Q_i$ with some added holding at
$\rho_i=\sigma^{-i+1}\rho\sigma^{i-1}$.  This is obviously a special
case of Example \ref{ex-wave-vertex'}, and we thus have
\begin{equation}\label{sticky}
\max\{\widetilde{\pi}\}\le c \min\{\widetilde{\pi}\},\qquad
c=\bigl(1-2n\delta/(n-1)\bigr)^{-1} .\end{equation}
Hence Proposition \ref{cor-wave} applies.  The second largest singular value
of $Q$ is  known to be $\sigma_1=1-1/(2n)$ (see, e.g., \cite{Dia2,FOW,SZ1}).
This yields an upper bound of order $n(n\log n +\log_+1/\eta)$
for the relative-sup merging time $T_\infty(\eta)$ of the sequence
$(K_i)_1^\infty$.  This result can be improved by using the
logarithmic Sobolev inequality technique of \cite{SZ4}, (\ref{sticky}) and
Lemma \ref{Dircomp}. The logarithmic Sobolev constant
$l(Q^2)$ of $Q^2$ is of order $1/n\log n$ (see \cite{DS-L}). This yields
a relative-sup merging time upper bound of order
$n((\log n)^2 +\log_+1/\eta)$.  This result holds also if we replace
the lazy random walk $Q$ above by its nonlazy version,
the usual ``transpose top with random.''

A total variation merging time estimate
of order $n(\log n +\log_+1/\eta)$ is obtained in \cite{SZ4} by using
Lemma \ref{lem-mMstab} together with the
modified logarithmic Sobolev inequality technique.
The crucial point is that the modified logarithmic Sobolev constant $l'(Q^2)$
of $Q^2$ is of order $1/n$ (see \cite{Goel,SZ4}). We do not know how to prove
this improved estimate for the nonlazy version of this example.
\end{example}


\printaddresses

\end{document}